\documentclass[fleqn]{mat01}
\usepackage{times,mathtimy,amssymb,latexsym}
\begin{document}

\setcounter{page}{331}
\firstpage{331}

\newtheorem{theore}{Theorem}
\renewcommand\thetheore{\arabic{theore}}
\newtheorem{theor}[theore]{\bf Theorem}
\newtheorem{propo}[theore]{\rm PROPOSITION}
\newtheorem{exam}[theore]{Example}

\def\exammp{\trivlist\item[\hskip\labelsep{\it Examples.}]}
\def\noot{\trivlist\item[\hskip\labelsep{\it Note.}]}

\title{Multipliers of $\pmb{A}_{\pmb{p}}\pmb{((0,}\,\pmb{\infty}\pmb{))}$ with order convolution}

\markboth{Savita Bhatnagar}{Multipliers of $A_{p}((0,\infty))$
with order convolution}

\author{SAVITA BHATNAGAR}

\address{Department of Mathematics, Panjab University,
Chandigarh~160~014, India\\
\noindent E-mail: bhsavita@pu.ac.in} 

\volume{115}

\mon{August}

\parts{3}

\pubyear{2005}

\Date{MS received 5 July 2004; revised 31 January 2005}

\begin{abstract}
The aim of this paper is to study the multipliers from $A_{r}(I)$
to $A_{p}(I), r \ne p$, where $I=(0,\infty)$ is the locally
compact topological semigroup with multiplication max and usual
topology and $A_{r}(I) = \{f \in L_{1}(I)\hbox{:}\ \hat{f} \in
L_{r}(\hat{I})\}$ with norm $|||f|||_{r} = \|f\|_{1} +
\|\hat{f}\|_{r}$.
\end{abstract}


\maketitle

\section{Introduction}

The algebra $A_{p}(G)$ of elements in $L_{1}(G)$ whose Fourier
transforms belong to $L_{p}(\hat{G})$ and the multipliers for
these algebras have been studied by various authors
\cite{1,7,8,9}. Let $I=(0,\infty)$ be the locally compact
idempotent commutative topological semigroup with the usual
topology and max multiplication and $\hat{I}$ be the maximal ideal
space of $L_{1}(I)$. The algebras $A_{p}(I)$ of elements in
$L_{1}(I)$ whose Gelfand transforms belong to $L_{p}(\hat{I})$ and
the $(A_{p}(I), A_{p}(I))$ multipliers for these algebras have
been studied by Kalra, Singh and Vasudeva \cite{4}. The purpose of
this note is to investigate the $(A_{r}(I), A_{p}(I))$ multipliers
for $1 \leq r, p \leq \infty$. Even in the group case, the study
of $(A_{r}(G), A_{p}(G))$ multipliers is not complete. Only some
partial results have been listed in \cite{2}. We have been able to
obtain a set of necessary conditions and another set of sufficient
conditions on a function $\varphi$ so that it defines a
multiplier. Some instructive examples which have bearing on the
above-said necessary and sufficient conditions have been provided.
The next section contains all the preliminary results which
we shall use throughout the paper. The last section contains the
results on the multipliers.

\section{Preliminaries}

Let $I=(0,\infty)$ be the locally compact semigroup with
multiplication max and usual topology. Let $M(I)$ denote the
Banach algebra of all finite regular Borel measures on $I$ under
the order convolution product denoted by $*$ and total variation
norm. Then the Banach space $L_{1}(I)$ of all measures in $M(I)$
which are absolutely continuous with respect to the Lebesgue
measure on $I$ becomes a commutative semisimple Banach algebra in
the inherited product $*$. More specifically, for $f, g \in
L_{1}(I)$,
\begin{equation*}
f * g(x) = f(x) \int_{0}^{x} g(y) {\rm d}y + g(x) \int_{0}^{x}
f(y) {\rm d}y\ \ \hbox{a.e.}
\end{equation*}
The maximal ideal space $\hat{I}$ of $L_{1}(I)$ can be identified
with the interval $(0,\infty]$ and the Gelfand transform $\hat{f}$
of an $f\in L_{1}(I)$ is the indefinite integral, i.e.,
\begin{equation*}
\hat{f}(x) = \int_{0}^{x} f(t) {\rm d}t.
\end{equation*}
For these and other results that may be used in the sequel, the
reader is referred to \cite{3,5}.

The algebras $A_{p}(I), 1\leq p \leq \infty$, consist of $f\in
L_{1}(I)$ such that $\hat{f}\in L_{p}(\hat{I})$. Clearly
$A_{\infty}(I) = L_{1}(I)$ and each $A_{p}(I)$ is an ideal in
$L_{1}(I)$. Define
\begin{equation*}
|||f|||_{p} = \|f\|_{1} + \|\hat{f}\|_{p}, \quad f\in A_{p}(I).
\end{equation*}
Then $|||\cdot |||_{p}$ is a norm on $A_{p}(I)$ and $A_{p}(I)$ is
a commutative Banach algebra with order convolution. Moreover,
$A_{p}(I)$ is a proper subset of $A_{r}(I),1\leq p < r\leq
\infty$. The maximal ideal space $\Delta(A_{p}(I))$ of $A_{p}(I)$
is homeomorphic to $I$, (note the contrast with the group case
$\Delta (A_{p}(G)) = \Delta (L_{1}(G)) \cong \hat{G}$, whereas in
the case of semigroup $I$ under consideration $\Delta(L_{1}(I)) =
\hat{I} \cong (0, \infty] \ne \Delta (A_{p}(I))$) and $A_{p}(I)$ is
a semisimple commutative Banach algebra which does not contain a
bounded approximate identity for $1\leq p < \infty$ but contains
an approximate identity. For the above results on $A_{p}(I)$, the
reader is referred to \cite{4}.

A mapping $T$ on a commutative Banach algebra $X$ to itself is
called a multiplier if $T(xy) = T(x)y$ for $x,y\in X$. If $X$ is
semisimple and $T\hbox{:}\ X \rightarrow X$ is a multiplier then
there exists a unique continuous and bounded function $\varphi$ on
$\Delta(X)$ such that $(\hat{T}x) = \varphi\hat{x}, x \in X$ and
$\|\varphi\|_{\infty} \leq \|T\|$ (p.~19 of \cite{6}). Due to
semisimplicity of $X, \hat{X}$ is also a commutative semisimple
Banach algebra under pointwise operations with norm
$\|\hat{x}\|=\|x\|$ and $\Delta (\hat{X})$ is homeomorphic to
$\Delta(X)$. Therefore, for semisimple Banach algebras we may
consider the multiplier $T$ to be an operator
$M_{\varphi}\hbox{:}\ \hat{X} \rightarrow \hat{X}, \hat{x} \in
\hat{X}$ defined by $M_{\varphi} (\hat{x}) = \varphi \hat{x},
\hat{x} \in \hat{X}$. Now $A_{p}(I)$ is semisimple, so
$\hat{A}_{p}(I) = \{\hat{f}\hbox{:}\ f\in A_{p}(I)\}$ is the
Banach algebra under pointwise operations and $|||\hat{f}||| =
|||f|||_{p} = \|f\|_{1} + \|\hat{f}\|_{p}$. A multiplier from
$A_{r}(I)$ to $A_{p}(I)$ may be defined (p.~67 of \cite{6}) as a
continuous function $\varphi$ on $\Delta (A_{r}(I)) = \Delta
(A_{p}(I)) \cong I$ such that $\varphi \hat{f} \in \hat{A}_{p}(I)$
whenever $\hat{f}\in \hat{A}_{r}(I)$.

Let $A$ denote the algebra of all complex-valued measurable
functions on $I$ under pointwise operations. For $\varphi \in A$,
let $M_{\varphi}$ denote the operator defined by
$M_{\varphi}(\hat{f}) = \varphi\hat{f}, \hat{f} \in
\hat{A}_{p}(I)$. The following theorem has been proved in
\cite{4}.

\begin{theor}[\!]
$M_{\varphi}$ is a bounded multiplier on $\hat{A}_{p}(I)$ to
itself iff
\begin{enumerate}
\renewcommand\labelenumi{\rm (\roman{enumi})}
\leftskip .35pc
\item $\varphi$ is bounded{\rm ,}

\item $\varphi$ is absolutely continuous on $[0,K]$ for each
$K>0${\rm ,}

\item $M_{\varphi'}$ is a bounded linear operator from
$\hat{A}_{p}(I)$ to $L_{1}(I)$.
\end{enumerate}
\end{theor}
In the next section, we provide a characterization of multipliers
from $A_{r}(I)$ to $A_{p}(I), r\ne p$. We give some intrinsic
conditions on $\varphi$ so that $\varphi$ defines a multiplier
from $A_{r}(I)$ to $A_{p}(I)$. Finally, we show that if $\varphi$
is a multiplier then $\varphi$ satisfies some growth conditions.

\section{The $\pmb{(}\pmb{A}_{\pmb{r}}\pmb{(I),}\ \pmb{A}_{\pmb{p}}\pmb{(I))}$ multipliers}

\begin{propo}$\left.\right.$\vspace{.5pc}

\noindent If $f\in A_{p}(I)$ and $|||f|||_{p} \leq 1${\rm ,} then
$f\in A_{r}(I), r >p$ and $|||f|||_{r} < 2$.
\end{propo}

\begin{proof}
If $f\in A_{p}(I)$ and $|||f|||_{p} = \|f\|_{1} + \|\hat{f}\|_{p}
\leq 1$, then $\|f\|_{1} < 1$ and $\|\hat{f}\|_{p} < 1$\break and
\begin{align*}
\|\hat{f}\|_{r}^{r} =\int_{0}^{\infty} |\hat{f}|^{r} &=
\int_{0}^{\infty} |\hat{f}|^{r-p} \cdot |\hat{f}|^{p}\\[.3pc]
&\leq \int_{0}^{\infty} |\hat{f}|^{p}\ \ \hbox{as}\ \ |\hat{f}(x)|
\leq \|f\|_{1} < 1 \ \ \hbox{and}\ \ r - p > 0\\[.3pc]
&= \|\hat{f}\|_{p}^{p} < 1.
\end{align*}
Thus $\|\hat{f}\|_{r} < 1$ and $|||f|||_{r} = \|f\|_{1} +
\|\hat{f}\|_{r} < 2$.\hfill $\Box$
\end{proof}

\begin{theor}[\!]
Let $M_{\varphi}\hbox{\rm :}\ \hat{A}_{r}(I) \rightarrow
\hat{A}_{p}(I), r > p$ be a bounded multiplier{\rm ,} then
\begin{enumerate}
\renewcommand\labelenumi{\rm (\roman{enumi})}
\leftskip .35pc
\item $\varphi$ is bounded{\rm ,}

\item $\varphi$ is absolutely continuous on $[0, K]$ for each $K
>0$ and
\begin{equation*}
\hskip -1.25pc \int_{0}^{x} |\varphi'(t)|{\rm d}t =
O(x^{1/r})\quad \hbox{as}\ \ x \rightarrow \infty,
\end{equation*}
\item $M_{\varphi'}$ is a bounded linear operator from
$\hat{A}_{r}(I)$ to $L_{1}(I)$.
\end{enumerate}
\end{theor}

\begin{proof}
Since $r > p, A_{p}(I) \subseteq A_{r}(I)$. Therefore if
$M_{\varphi}\hbox{:}\ \hat{A}_{r}(I) \rightarrow \hat{A}_{p}(I)$
is a bounded multiplier then $M_{\varphi}$ defines a bounded
multiplier from $\hat{A}_{r}(I)$ to itself. The theorem now
follows from Theorems~15 and 18 of \cite{4}. It should be noted
that $\|\varphi\|_{\infty} \leq 2 \|M_{\varphi}\|$. To see this,
let
\begin{equation*}
K_{x} = \mathop{\sup}\limits_{|||\hat{f}|||_{p} \leq 1}
|\hat{f}(x)|,\quad x \in (0, \infty).
\end{equation*}
Then
\begin{equation*}
K_{x} \leq \mathop{\sup}\limits_{|||\hat{f}|||_{r}\leq 2}
|\hat{f}(x)| = 2 \mathop{\sup}\limits_{|||\hat{f}/2|||_{r}\leq 1}
|(\hat{f}/2)(x)| = 2 \mathop{\sup}\limits_{|||\hat{f}|||_{r}\leq
1} |\hat{f}(x)|,
\end{equation*}
using Proposition~2.

Now,
\begin{equation*}
|\varphi(x) \hat{f}(x)| = |M_{\varphi}(\hat{f})| \leq
K_{x}\|M_{\varphi}\|\ |||\hat{f}|||_{r}.
\end{equation*}
Therefore, in particular,
\begin{align*}
|\varphi(x)| \leq \mathop{\inf}\limits_{|||\hat{f}|||_{r}\leq 1}
\frac{K_{x}\|M_{\varphi}\|}{|\hat{f}(x)|} =
\frac{K_{x}\|M_{\varphi}\|}{\mathop{\sup}_{|||\hat{f}|||_{r}\leq
1} |\hat{f}(x)|} \leq \frac{2K_{x}\|M_{\varphi}\|}{K_{x}} =
2\|M_{\varphi}\|.
\end{align*}

$\left.\right.$\vspace{-1.5pc}

\hfill $\Box$
\end{proof}
The next theorem gives sufficient conditions on $\varphi$ so that
$\varphi$ defines a multiplier from $A_{r}(I)$ to
$A_{p}(I),r>p$.\pagebreak

\begin{theor}[\!]
Let $r > p$ and $1/v=1/p-1/r >0$. If
\begin{enumerate}
\renewcommand\labelenumi{\rm (\roman{enumi})}
\leftskip .35pc
\item $\varphi \in L_{v}(I)${\rm ,}

\item $\varphi$ is absolutely continuous on every interval $[0,K],
K > 0${\rm ,}

\item $M_{\varphi'}$ is a bounded multiplier on $\hat{A}_{r}(I)$
to $L_{1}(I)${\rm ,}\vspace{-.5pc}
\end{enumerate}
then $M_{\varphi}\hbox{\rm :}\ \hat{A}_{r}(I) \rightarrow
\hat{A}_{p}(I)$ is a bounded multiplier.
\end{theor}

\begin{proof}
Suppose $\varphi$ satisfies (i), (ii) and (iii). By (i),
$\lim_{x\rightarrow \infty} \varphi(x) = 0$ so there exists $K$
such that for all $x > K, |\varphi(x)|< 1$. Using (ii) we get
$\varphi$ to be bounded on $[0, K]$. Thus $\varphi$ is bounded.
Since $\varphi$ and $\hat{f}$ are absolutely continuous on $[0,K]$
for each $K> 0$ so is $\varphi \hat{f}$. Thus the derivative of
$\varphi\hat{f}, (\varphi\hat{f})'$ exists a.e. on $I$ and
$(\varphi\hat{f})' = \varphi'\hat{f} + \varphi f\in L_{1}(I)$ in
view of (iii) and boundedness of $\varphi$. Moreover,
\begin{align*}
\|(\varphi\hat{f})'\|_{1} &= \|\varphi'\hat{f} + \varphi f\|_{1}\\[.3pc]
&\leq \|M_{\varphi'}\|\ |||\hat{f}|||_{r} +
\|\varphi\|_{\infty}\|f\|_{1}.
\end{align*}
Since $\varphi \in L_{v}(I)$ and $\hat{f}\in L_{r}(I)$, it follows
that $\varphi\hat{f} \in L_{p}(I)$ and
\begin{equation*}
\|\varphi\hat{f}\|_{p}\leq \|\varphi\|_{v}\|\hat{f}\|_{r}.
\end{equation*}
Thus $\varphi\hat{f}\in \hat{A}_{p}(I)$ and
\begin{align*}
|||\varphi\hat{f}|||_{p} &= \|(\varphi\hat{f})'\|_{1} + \|\varphi
\hat{f}\|_{p}\\[.3pc]
&\leq (\|M_{\varphi'}\| + \|\varphi\|_{\infty} + \|\varphi\|_{v})
|||f|||_{r}.
\end{align*}
$M_{\varphi}$ is a bounded multiplier on $\hat{A}_{r}(I)$ to
$\hat{A}_{p}(I)$ and
\begin{equation*}
\|M_{\varphi}\| \leq \|M_{\varphi'}\| + \|\varphi\|_{\infty} +
\|\varphi\|_{v}.
\end{equation*}
This completes the proof.\hfill $\Box$
\end{proof}

Note that (iii) may be replaced by $\varphi'|_{(a,\infty)} \in
L_{q}(a, \infty)$ for some $a > 0$ and some $q, 1 \leq q \leq r'$.
In that case
\begin{equation*}
\hat{f} \in L_{r} (a, \infty) \subseteq L_{q'} (a, \infty)
\end{equation*}
and
\begin{equation*}
\|\hat{f}|_{(a,\infty)}\|_{q'}^{q'} \leq
\|\hat{f}|_{(a,\infty)}\|_{r}^{r} \cdot
\|\hat{f}\|_{\infty}^{q'-r} \leq \|\hat{f}|_{(a,\infty)}\|_{r}^{r}
\cdot \|f\|_{1}^{q'-r}
\end{equation*}
so that $\|M_{\varphi}\| \leq \|\varphi'|_{[0,a]}\|_{1} +
\|\varphi'|_{(a,\infty)}\|_{q} + \|\varphi\|_{\infty} +
\|\varphi\|_{v}$.

\begin{exam}$\left.\right.$
{\rm
\begin{enumerate}
\renewcommand\labelenumi{\rm (\roman{enumi})}
\leftskip .2pc
\item For $r> p \geq 1, 1/v = 1/p-1/r$. Take $\varphi(x) = 1, x
\in I$.
Then $\varphi$ is bounded and continuous on $I$ but $\varphi
\notin L_{v}(I), \varphi' = 0 \in L_{1}(I)$. So (ii) and (iii)
hold but (i) does not hold in Theorem~4. For
\begin{align*}
\hskip -1.25pc f(x) &= \begin{cases}
1, &0 < x < 1,\\[.3pc]
-\alpha x^{-\alpha-1}, &x \geq 1,
\end{cases}
\end{align*}
\begin{align*}
\hskip -1.25pc \hat{f}(x) &= \begin{cases}
x, &0 < x < 1,\\[.3pc]
x^{-\alpha}, &x \geq 1.
\end{cases}
\end{align*}
For $0 < 1/r < \alpha < 1/p$ ($\alpha$ exists as $r > p$),
$\hat{f} \in L_{r}(I)\backslash L_{p}(I)$, so $\varphi\hat{f} =
\hat{f} \notin L_{p} (I)$. Thus $\varphi$ is not a multiplier on
$\hat{A}_{r}(I)$ to $\hat{A}_{p}(I)$.

\item If $r > p \geq 1, 1/v = 1/p-1/r >0$ and $\epsilon >0$, take
\begin{equation*}
\hskip -1.25pc \varphi(x) =
\begin{cases}
1, &0 < x \leq 1, \\[.3pc]
x^{-\frac{1}{v}-\epsilon}, &x > 1.
\end{cases}
\end{equation*}
Then $\varphi$ is continuous on $I, \varphi \in L_{v}(I)$,
$\varphi$ is absolutely continuous on $[0,K]$ for all $K>0$.
\begin{align*}
\hskip -1.25pc \varphi'(x) &= \begin{cases}
0, &0 < x < 1, \\[.3pc]
\left( -\frac{1}{v} - \epsilon \right)x^{-\frac{1}{v}-\epsilon-1},
&x > 1,
\end{cases}\\[.3pc]
\hskip -1.25pc &= \begin{cases}
0, &0 < x < 1, \\[.3pc]
\left( -\frac{1}{v} - \epsilon \right)
x^{-\frac{1}{r'}-\frac{1}{p}-\epsilon}, &x > 1,
\end{cases}
\end{align*}
is in $L_{r'}(1, \infty)$ so that conditions of Theorem~4 are
satisfied. Hence $M_{\varphi}$ is a multiplier from
$\hat{A}_{r}(I)$ to $\hat{A}_{p}(I)$.
\end{enumerate}}
\end{exam}

In the next two theorems, we discuss the case $r < p$.

\begin{theor}[\!]
If
\begin{enumerate}
\renewcommand\labelenumi{\rm (\roman{enumi})}
\leftskip .35pc
\item $\varphi$ is bounded on $(0,\infty)$ or $\varphi \in
L_{p}(I)${\rm ,}

\item $\varphi$ is absolutely continuous on $[0,K]$ for each $K >
0${\rm ,}

\item $M_{\varphi'}$ is a bounded multiplier on $\hat{A}_{r}(I)$
to $L_{1}(I)${\rm ,} then $M_{\varphi}\hbox{\rm :}\ \hat{A}_{r}(I)
\rightarrow \hat{A}_{p}(I), r < p$ is a bounded multiplier.
\end{enumerate}
\end{theor}

\begin{proof}
If $\varphi$ is bounded on $I$ and (ii) and (iii) hold, then
$M_{\varphi}\hbox{:}\ \hat{A}_{r}(I) \rightarrow \hat{A}_{r}(I)$
is a bounded multiplier by Theorem~1. So $M_{\varphi}\hbox{:}\
\hat{A}_{r}(I) \rightarrow \hat{A}_{p}(I)$ is a bounded multiplier
as $\hat{A}_{r}(I) \subseteq \hat{A}_{p}(I)$ for $r<p$. Moreover,
\begin{align*}
\|M_{\varphi}\| &= \mathop{\sup}\limits_{|||\hat{f}|||_{r}= 1}
(|||\varphi\hat{f}|||_{p})\\[.3pc]
&= \mathop{\sup}\limits_{|||\hat{f}|||_{r}= 1}
(\|\varphi\hat{f}\|_{p} +
\|\varphi'\hat{f} + \varphi f \|_{1})\\[.3pc]
&\leq \mathop{\sup}\limits_{|||\hat{f}|||_{r}= 1}
(\|\varphi\|_{\infty} \|\hat{f}\|_{p} + \|M_{\varphi'}\|\
|||\hat{f}|||_{r} + \|\varphi\|_{\infty} \|f\|_{1})\\[.3pc]
&\leq 2 \|\varphi\|_{\infty} + \|M_{\varphi'}\|\ \hbox{as}\
\|f\|_{1} < 1\ \,\hbox{and}\ \,\|\hat{f}\|_{p} < 1 \ \hbox{using
Proposition~2.}
\end{align*}
Next, if $\varphi \in L_{p}(I)$, then as in the proof of
Theorem~4, $\varphi$ is bounded on $I$ and hence
$M_{\varphi}\hbox{:}\ \hat{A}_{r}(I) \rightarrow \hat{A}_{r}(I)$
is a bounded multiplier (using Theorem~1).\hfill $\Box$
\end{proof}

\begin{theor}[\!]
If $M_{\varphi}$ is a bounded multiplier on $\hat{A}_{r}(I)$ to
$\hat{A}_{p}(I), r < p$ then $\varphi$ is continuous on $I,
\varphi$ is absolutely continuous on $[0,K]$ for $K>0, \varphi$ is
locally in $L_{p}(I)$ and
\begin{align*}
\|\varphi|_{[0,x]}\|_{p} &= o\left(x^{\frac{1}{r} + \epsilon}
\right),\\[.3pc]
\|\varphi'|_{[0,x]}\|_{1} &= o\left(x^{\frac{1}{r} + \epsilon}
\right)\quad \hbox{for all}\ \ \epsilon >0.
\end{align*}
\end{theor}

\begin{proof}
The proof is exactly similar to that of the proof of theorem~15 of 
\cite{4}. Let $M_{\varphi}\hbox{:}\ \hat{A}_{r}(I) \rightarrow
\hat{A}_{p}(I), r < p$ be a bounded multiplier, then $\varphi$ is
continuous and
\begin{equation*}
|||\varphi\hat{f}|||_{p} = |||M_{\varphi}(\hat{f})|||_{p} \leq
\|M_{\varphi}\|\ |||\hat{f}|||_{r},\quad f \in A_{r}(I).
\end{equation*}
Take the function
\begin{equation*}
f(x) = \begin{cases}
1/\alpha, &0 < x < \alpha, \\[.3pc]
0, &\alpha \leq x < \beta,\\[.3pc]
1/(\beta-\gamma), &\beta \leq x < \gamma,\\[.3pc]
0, &x \geq \gamma.
\end{cases}
\end{equation*}
As in the proof of theorem~15 of \cite{4}, $\varphi\hat{f} \in
\hat{A}_{p}(I)$ implies $\varphi$ is continuous on $(0,\infty),
\varphi'$ exists a.e. on $(0,\infty)$ and
\begin{equation*}
\|\varphi \hat{f}\|_{p} + \|\varphi'\hat{f} + \varphi f\|_{1} =
|||\varphi \hat{f}|||_{p} \leq \|M_{\varphi}\|\ |||\hat{f}|||_{r}.
\end{equation*}
This implies that
\begin{equation*}
\|\varphi|_{[\alpha,\beta]}\|_{p} +
\|\varphi'|_{[\alpha,\beta]}\|_{1} \leq \|M_{\varphi}\| \left[ 2 +
\left( \frac{\gamma + r (\beta - \alpha)}{r + 1} \right)^{1/r}
\right].
\end{equation*}
Fix $\beta$ and vary $\gamma > \beta$. Then
\begin{equation*}
\mathop{\rm sup}\limits_{0<\alpha < \beta}
(\|\varphi|_{[\alpha,\beta]}\|_{p} + \| \varphi'|_{[\alpha,
\beta]}\|_{1}) \leq \|M_{\varphi}\| (2 + \beta^{1/r}).
\end{equation*}
So $\|\varphi|_{[0,\beta]}\|_{p} = o (\beta^{\frac{1}{r}+
\epsilon})$ and $\|\varphi'|_{[0,\beta]}\|_{1} = o
(\beta^{\frac{1}{r}+ \epsilon})$ for all $\epsilon > 0$. Thus
$\varphi$ is absolutely continuous on $[0,\beta]$ for $\beta >0$.
This completes the proof.\hfill $\Box$
\end{proof}

\begin{exam}{\rm
For $\epsilon > 0$ and $r < p$, consider
\begin{equation*}
\varphi(x) =
\begin{cases}
1, &0 < x < 1, \\[.3pc]
x^{\frac{1}{r}-\frac{1}{p}+\epsilon}, &x\geq 1.
\end{cases}
\end{equation*}
Then
\begin{equation*}
\varphi'(x) =
\begin{cases}
0, &0 < x < 1, \\[.3pc]
\left( \frac{1}{r} - \frac{1}{p} + \epsilon \right)
x^{\frac{1}{r}-\frac{1}{p}+\epsilon-1}, &x\geq 1,
\end{cases}
\end{equation*}
$\varphi$ is continuous on $I$, $\varphi'$ exists a.e. on $I$ and
\begin{equation*}
\|\varphi'|_{[0,\beta]}\|_{1} = \beta^{\frac{1}{r} - \frac{1}{p} +
\epsilon} - 1 = o \left(\beta^{\frac{1}{r}+ 2\epsilon}\right)
\end{equation*}
and
\begin{equation*}
\|\varphi|_{[0,\beta]}\|_{p} = \left(1 +
\frac{\beta^{p/r+p\epsilon}-1}{p/r+p\epsilon} \right)^{1/p} = o
\left(\beta^{\frac{1}{r}+ 2\epsilon}\right) \quad \hbox{for}\ \
\epsilon > 0.
\end{equation*}
Thus $\varphi$ satisfies the necessary conditions of Theorem~7 for
being a multiplier from $\hat{A}_{r}(I)$ to $\hat{A}_{p}(I), r <
p$. But $\varphi$ is not bounded and $\varphi \notin L_{p}(I)$.
Also $\varphi' \in L_{r'}(I)$ for $0 < \epsilon < 1/p$. It follows
that $M_{\varphi'}(I)$ is a bounded multiplier from
$\hat{A}_{r}(I)$ to $L_{1}(I)$ and (iii) holds in Theorem~6.

For
\begin{align*}
f(x) &=
\begin{cases}
1, &0 < x < 1, \\[.3pc]
\left( - \frac{1}{r} - \epsilon\right) x^{-\frac{1}{r} - \epsilon -1},
&x \geq 1,
\end{cases} \\[.3pc]
\hat{f}(x) &=
\begin{cases}
x, &0 < x < 1, \\[.3pc]
x^{-\frac{1}{r}-\epsilon}, &x\geq 1,
\end{cases}
\end{align*}
$f\in A_{r}(I)$, but
\begin{equation*}
\varphi \hat{f}(x) = \begin{cases}
x, &0 < x < 1, \\[.3pc]
x^{-\frac{1}{p}}, &x \geq 1,
\end{cases}
\end{equation*}
is not in $L_{p}(I)$. So $M_{\varphi}$ is not a multiplier from
$\hat{A}_{r}(I)$ to $\hat{A}_{p}(I)$.}
\end{exam}

\begin{noot}
All the multipliers from $\hat{A}_{r}(I)$ to $\hat{A}_{p}(I), r <
p$, will not be given by continuous bounded functions $\varphi$
otherwise as in the proof of theorem~1 bounded
multipliers from $\hat{A}_{r}(I)$ to $\hat{A}_{p}(I), r < p$ will
only be the multipliers from $\hat{A}_{r}(I)$ to $\hat{A}_{r}(I)
\subset \hat{A}_{p}(I)$.
\end{noot}

\section*{Acknowledgement}

The author is grateful to Prof.~H~L~Vasudeva for useful
discussions.

\end{document}